\documentclass[preprint,12pt]{article}
\usepackage{}
\usepackage{amssymb}
\usepackage{mathrsfs}
\usepackage{amsfonts}
\usepackage{graphicx}
\usepackage{amsmath}
\usepackage{slashbox}
\usepackage{cite}
\usepackage{subfigure}
\newtheorem{thm}{Theorem}[subsection]

\newtheorem{lem}[thm]{Lemma}
\newtheorem{prop}[thm]{Proposition}

\def\QEDclosed{\mbox{\rule[0pt]{1.3ex}{1.3ex}}}
\numberwithin{equation}{section}
\numberwithin{thm}{section}

\begin{document}
\title{Effective Computation of Stochastic Protein Kinetic Equation by Reducing Stiffness via Variable Transformation}
       \author{
        Lijin Wang\footnotemark[1]\\
       {\small \footnotemark[1] School of Mathematical Sciences, University of Chinese Academy of Sciences,}\\
    {\small Beijing 100049, P.R.China }}
       \maketitle

        \footnotetext{ljwang@ucas.ac.cn}

       \begin{abstract}
          {\rm\small The stochastic protein kinetic equations can be stiff for certain parameters, which makes their numerical simulation rely on very small time step sizes, resulting in large computational cost and accumulated round-off errors. For such situation, we provide a method of reducing stiffness of the stochastic protein kinetic equation by means of a kind of variable transformation. Theoretical and numerical analysis show effectiveness of this method. Its generalization to a more general class of stochastic differential equation models is also discussed.}\\

\textbf{AMS subject classification: } {\rm\small 65C20, 65C30.}\\

\textbf{Key Words: }{\rm\small} Numerical methods for SDEs; Stochastic protein kinetic equation; Stiffness; Midpoint rule.
\end{abstract}

\section{Introduction}
\label{intro}
Consider the following stochastic differential equation describing the kinetics of the proportion $x$ of one of two possible forms of certain proteins
\begin{equation}\label{1.1}
dx=(1-x+\lambda x(1-x))dt+\sigma x(1-x)\circ dW(t),\quad x(0)=x_0,
\end{equation}
where $\lambda$ is interaction coefficient of the two proteins, $\sigma$ is the amplitude of the random Gaussian perturbation, and $W(t)$ is the standard Brownian motion. The small circle $\circ$ before $dW(t)$ denotes the stochastic integral of Stratonovich sense (\cite{kandp}). There is no explicit solution to this equation, wherefore numerical computations simulating the propagation of $x$ is needed. One of the most well-known numerical methods for solving stochastic differential equations (SDEs) is the Euler-Maruyama method, which is, however, only consistent to SDEs of It\^{o} type. For SDEs of Stratonovich type, the consistent method is the midpoint rule, which, when applied to (\ref{1.1}), takes the form:
\begin{equation}\label{1.2}
\begin{split}
x_{n+1}&=x_n+h(1+(\lambda-1)\frac{x_n+x_{n+1}}{2}-\lambda\frac{(x_n+x_{n+1})^2}{4})\\
&+\Delta W_n(\sigma\frac{x_n+x_{n+1}}{2}-\sigma\frac{(x_n+x_{n+1})^2}{4}),
\end{split}
\end{equation}
where $h=t_{n+1}-t_n$ is the uniform time step size, and $\Delta W_n=W(t_{n+1})-W(t_n)$ obeying the Gaussian distribution $\mathcal{N}(0,h)$ and independent for different $n$. The mean-square convergence of this method can be proved (\cite{kandp,mil}).

\begin{figure*}
\subfigure[$h=0.13$]{\label{f1.1}\includegraphics[width=0.3\textwidth]{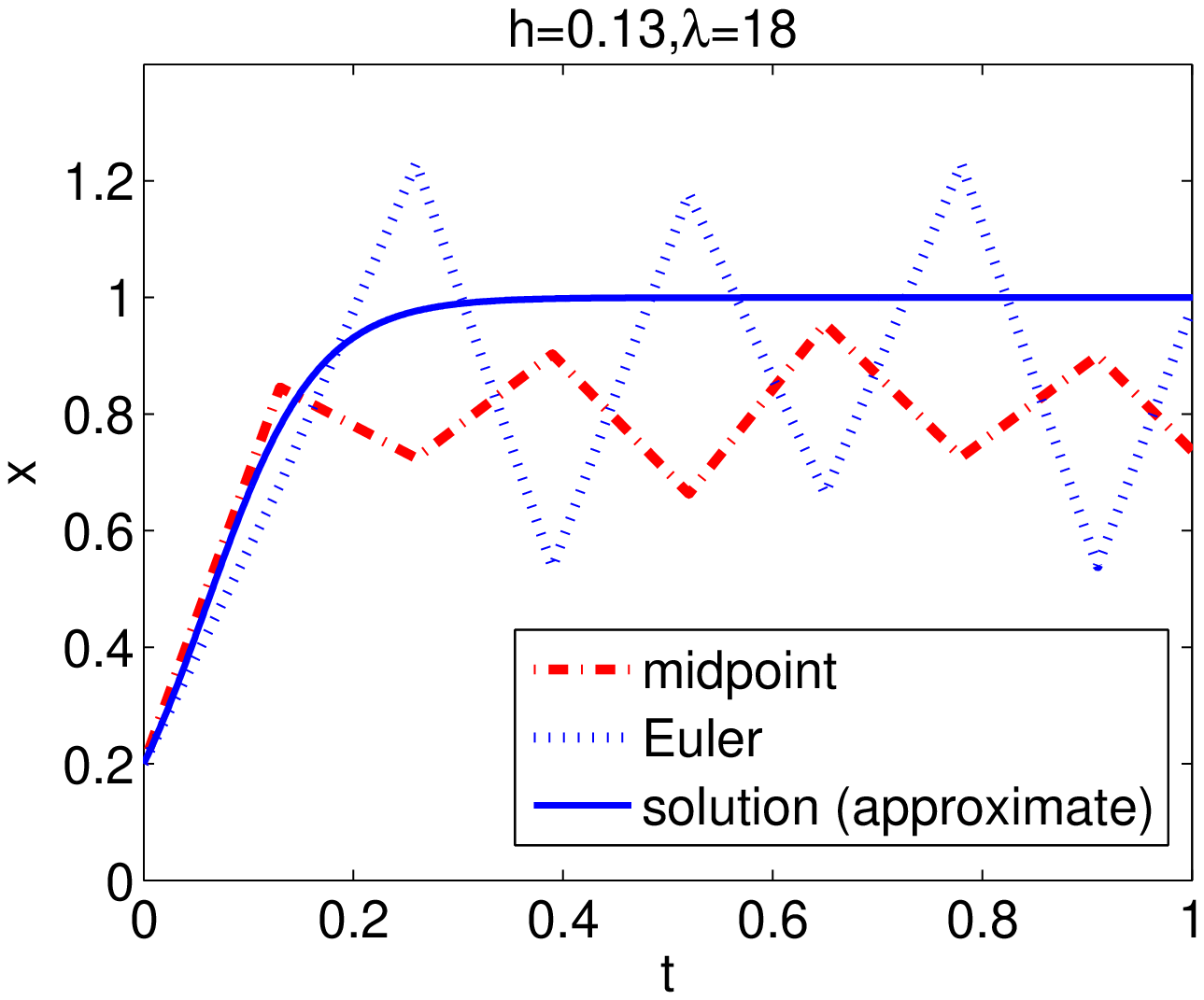}}
\hspace{1mm}
\subfigure[$h=0.1$]{\label{f1.2}\includegraphics[width=0.3\textwidth]{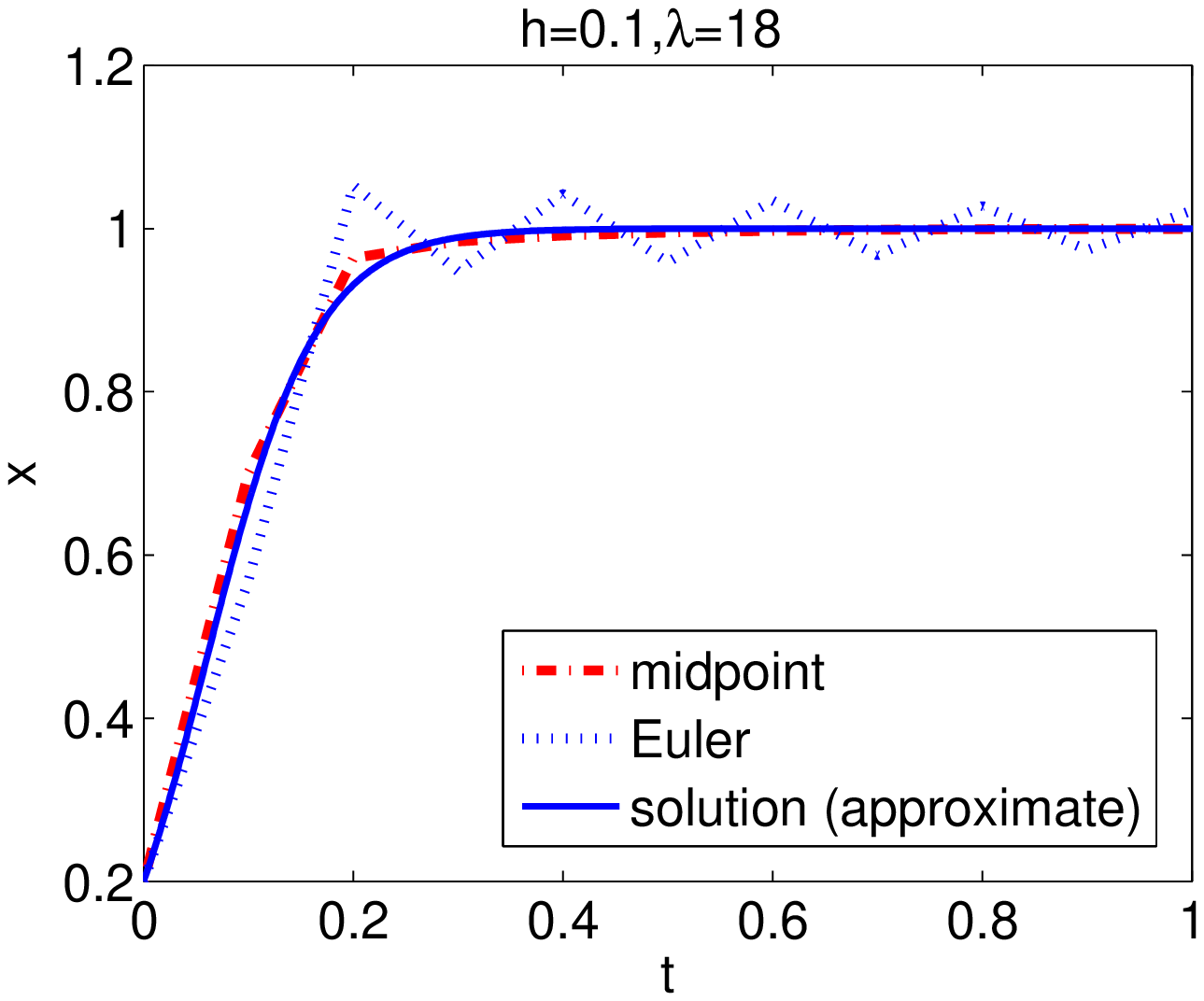}}
\hspace{1mm}
\subfigure[$h=0.01$]{\label{f1.3}\includegraphics[width=0.3\textwidth]{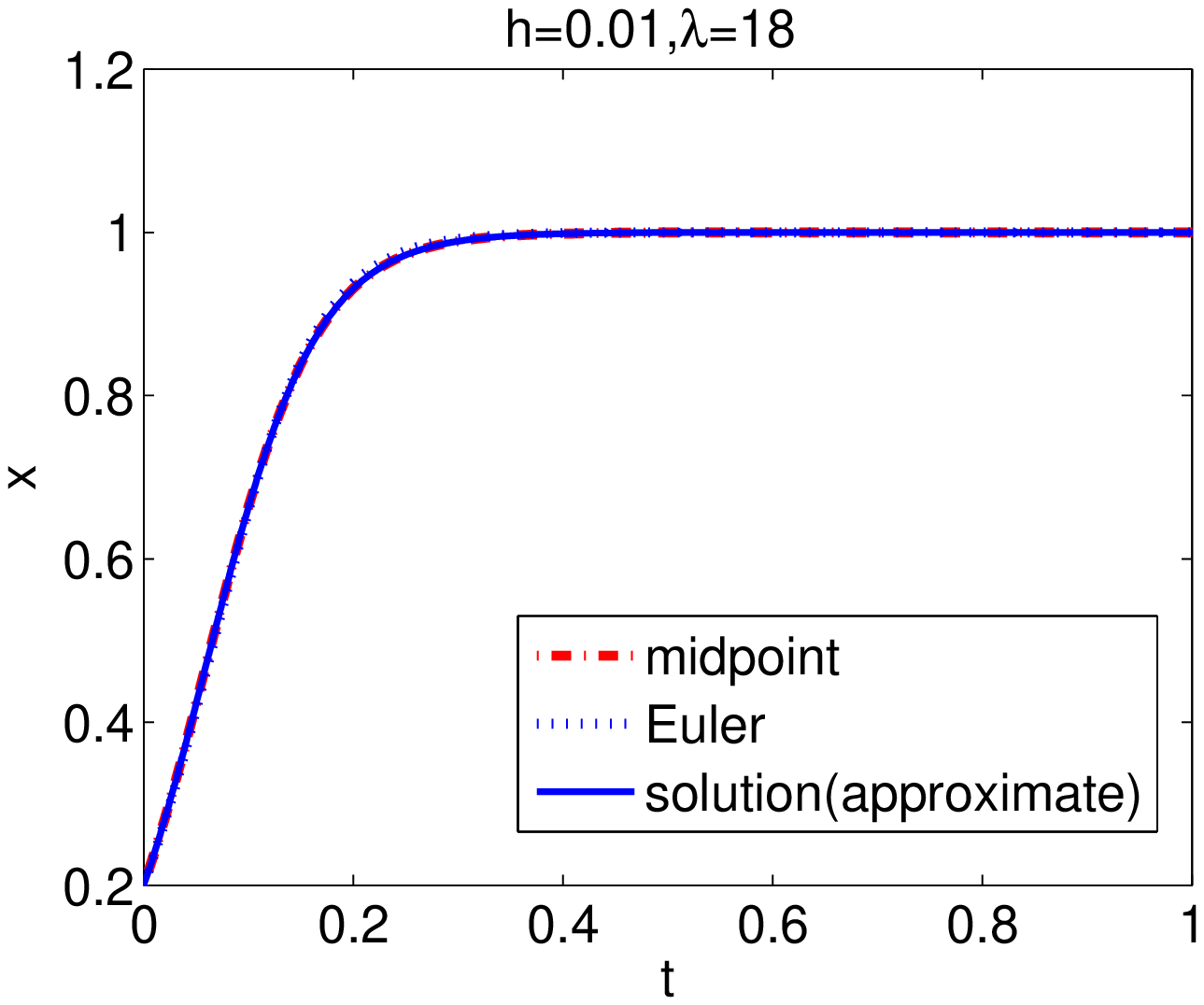}}
\caption{Illustration of stiffness and numerical stability}
\label{f1}       
\end{figure*}

The deterministic midpoint rule is an A-stable method suitable for dealing with stiff equations. For example, for the deterministic version of the equation (\ref{1.1}), i.e., the equation (\ref{1.1}) with $\sigma=0$
\begin{equation}\label{1.3}
\dot{x}=1-x+\lambda x(1-x),\quad x(0)=x_0,
\end{equation}
$x=1$ is an asymptotically stable solution as $\lambda>-1$. With $y=1-x$ we get the linearized equation of (\ref{1.3})
\begin{equation}\label{1.4}
\dot{y}=(-\lambda-1)y,
\end{equation}
which can be very stiff as $|\lambda|$ is large, e.g., $\lambda=18$. Then in this case, the original non-linear equation (\ref{1.3}) is also stiff, meaning that its numerical simulation may need the choice of very small step size $h$. To illustrate this, we use the Euler method
\begin{equation}
x_{n+1}=x_n+h(1+(\lambda-1)x_n-\lambda x_n^2)
\end{equation}
and the midpoint rule
\begin{equation}
x_{n+1}=x_n+h(1+(\lambda-1)\frac{x_n+x_{n+1}}{2}-\lambda\frac{(x_n+x_{n+1})^2}{4})
\end{equation}
to solve (\ref{1.3}) numerically, and observe the effect in Fig. \ref{f1}.

In Fig. \ref{f1.1}, $h=0.13$, we see that both the Euler method (blue dotted) and the midpoint rule (red dash-dotted) produce significant oscillating simulation error, while in Fig. \ref{f1.2}, where $h=0.1$, the midpoint rule creates a reasonable numerical solution remaining close to the true solution. Note that there is no explicit true solution for (\ref{1.3}), we just use the midpoint rule with a tiny step size $h=0.001$ to approximate the true solution (blue solid). However, if the step size is small enough, e.g. $h=0.01$ as used in Fig. \ref{f1.3}, both of the two methods give very accurate simulations which coincide visually with the true solution. In fact, as $|\lambda|$ gets larger and larger, the needed $h$ for effective computation becomes smaller and smaller. This illustrates that the simulation of the equation (\ref{1.3}) is sensitive to the time step size $h$, due to stiffness of the equation. Meanwhile, it can be seen that the midpoint rule is more stable than the Euler method.

For the stochastic protein kinetic equation (\ref{1.1}), $x=1$ is also a stochastically asymptotically stable solution if $\lambda>-1$. This can be seen from its linearized stochastic differential equation with $z=1-x$ (\cite{kandp})
\begin{equation}\label{1.7}
dz=(-\lambda-1)zdt-\sigma z\circ dW(t),
\end{equation}
the solution of which is $z_t=z_0 \exp((-\lambda-1)t-\sigma W(t))$. Thus the Lyapunov exponent $\mu$ of (\ref{1.7}) is (\cite{kandp,mao})
\begin{equation}\label{1.8}
\mu=\limsup_{t\rightarrow \infty}\frac{\ln |z|}{t}=-\lambda-1
\end{equation}
due to $\lim_{t\rightarrow\infty} \frac{W_t}{t}=0$ with probability 1.

At the same time, similar to its deterministic counterpart (\ref{1.4}), the equation (\ref{1.7}) can also be very sensitive to its simulation time step size $h$ as $|\lambda|$ is large, e.g. $\lambda=18$. In other words, we say (\ref{1.7}) is stiff for large $|\lambda|$. In this case, the original non-linear equation (\ref{1.1}) is also stiff (\cite{kandp}). This can be seen by numerical experiments, which we show in section 3, where we see that although the stochastic midpoint rule (\ref{1.2}) is applied, cautious choice of small time step size $h$ is still needed, which increases computational cost and accumulates round-off errors, especially for large time intervals.

Therefore it is meaningful to investigate effective measures to reduce the stiffness of the stochastic differential equations, such as the equation (\ref{1.1}), for efficient numerical simulation of such equations.
\section{Methods}
\label{sec:1}
For the purpose mentioned above, we employ the technique of variable transformation.
\begin{lem}\label{lem2.1} For $\lambda>-1$, the linearized stochastic differential equation (\ref{1.7}) is equivalent to the following linear SDE
\begin{equation}\label{2.1}
dZ=-Zdt-\frac{\sigma}{1+\lambda}Z\circ dW(t)
\end{equation}
via the variable transformation
\begin{equation}\label{2.2}
Z=z^{\frac{1}{1+\lambda}}.
\end{equation}
\end{lem}

{\bf Proof.}
Since the ordinary differential chain rule holds for SDEs of Stratonovich type, a direct calculation yields
\begin{equation}
\begin{split}
dZ&=\frac{1}{1+\lambda}z^{\frac{1}{1+\lambda}-1}dz\\
&=\frac{1}{1+\lambda}\frac{Z}{z}[-(\lambda+1)zdt-\sigma z\circ dW(t)]\\
&=-Zdt-\frac{\sigma}{1+\lambda}Z\circ dW(t).
\end{split}
\end{equation}
\hfill \QEDclosed

Note that, for large $|\lambda|$, the equation (\ref{1.7}) is much more stiff than its equivalent equation (\ref{2.1}), since the Lyapunov exponent for (\ref{2.1}) is $-1$, obtained in the same way as in (\ref{1.8}). We thus find a way of reducing the stiffness of (\ref{1.7}) via the variable transformation (\ref{2.2}). Here, $\lambda>-1$ is to guarantee that $z=0$ is the stochastically asymptotically stable solution of (\ref{1.7}), and that the transformation (\ref{2.2}) is meaningful at $z=0$.

Applying the same transformation to the original non-linear stochastic protein kinetic equation (\ref{1.1}), we obtain the following result.
\begin{prop}\label{prop2.2} For $\lambda>-1$, the stochastic protein kinetic equation (\ref{1.1}) is equivalent to the following SDE
\begin{equation}\label{2.4}
\begin{split}
dX&=(-X+\frac{\lambda}{1+\lambda}X^{2+\lambda})dt-(\frac{\sigma}{1+\lambda}X-\frac{\sigma}{1+\lambda}X^{2+\lambda})\circ dW(t),\\
X(0)&=(1-x_0)^{\frac{1}{1+\lambda}}
\end{split}
\end{equation}
via the variable transformation
\begin{equation}\label{2.5}
X=(1-x)^{\frac{1}{1+\lambda}}.
\end{equation}
The linearized equation of (\ref{2.4}) at its stochastically asymptotically stable solution $X=0$ is the equation (\ref{2.1}).
\end{prop}

{\bf Proof.} The inverse transformation of (\ref{2.5}) is
\begin{equation}\label{2.6}
x=1-X^{1+\lambda}.
\end{equation}
Using the differential chain rule together with (\ref{2.5})-(\ref{2.6}) we obtain
\begin{equation}
\begin{split}
dX&=\frac{-1}{1+\lambda}(1-x)^{\frac{1}{1+\lambda}-1}dx\\
&=\frac{-1}{1+\lambda}\frac{X}{1-x}[(1-x+\lambda x(1-x))dt+\sigma x(1-x)\circ dW(t)]\\
&=\frac{-1}{1+\lambda}[(X+\lambda xX)dt+\sigma xX\circ dW(t)]\\
&=(-X+\frac{\lambda}{1+\lambda}X^{2+\lambda})dt-(\frac{\sigma}{1+\lambda}X-\frac{\sigma}{1+\lambda}X^{2+\lambda})\circ dW(t).
\end{split}
\end{equation}
The condition $\lambda>-1$ can ensures that the transformation (\ref{2.5}) and its inverse (\ref{2.6}) are meaningful at $x=1$ and $X=0$, respectively. Meanwhile, it makes $X=0$ a stochastically asymptotically stable solution of (\ref{2.4}). The linearized equation of (\ref{2.4}) at $X=0$ is
\begin{equation}
\begin{split}
dZ&=(-1+\frac{\lambda(2+\lambda)}{1+\lambda}X^{1+\lambda}\left|_{X=0}\right.)Zdt\\
&-(\frac{\sigma}{1+\lambda}-
\frac{\sigma(2+\lambda)}{1+\lambda}X^{1+\lambda}\left|_{X=0}\right.)\circ dW(t)\\
&=-Zdt-\frac{\sigma}{1+\lambda}Z\circ dW(t),
\end{split}
\end{equation}
which is just the equation (\ref{2.1}). \hfill \QEDclosed \\

From Proposition \ref{prop2.2}, it is clear that the equation (\ref{2.4}) is much less stiff than the original equation (\ref{1.1}), since the linearized equation (\ref{2.1}) of (\ref{2.4}) is much less stiff than the linearized equation (\ref{1.7}) of (\ref{1.1}) for $|\lambda|$ large (\cite{kandp}). Thus, for the simulation of the stochastic protein kinetic equation, we can firstly apply the stochastic midpoint rule to the transformed equation (\ref{2.4}) to get $\{X_n\}$, and then use the inverse transform (\ref{2.6}) to get back to $\{x_n\}$.

In the following section, we perform numerical tests to illustrate the difference the transformation (\ref{2.5}) makes, that is, different degree of stiffness of (\ref{1.1}) and (\ref{2.4}), reflected in the dependence on time step size $h$ in the numerical simulations of them using the stochastic midpoint rule.
\section{Results}
\begin{figure*}
\subfigure[$h=0.12,\quad\lambda=18$]{\label{f2.1}\includegraphics[width=0.5\textwidth]{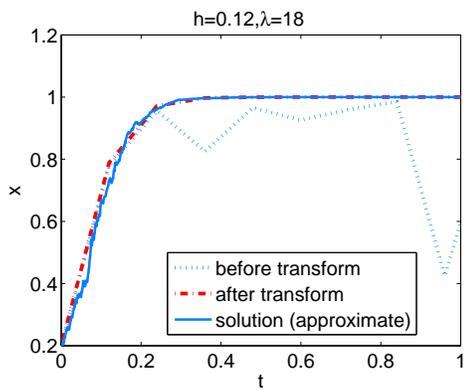}}
\hspace{1mm}
\subfigure[$h=0.01,\quad\lambda=18$]{\label{f2.2}\includegraphics[width=0.5\textwidth]{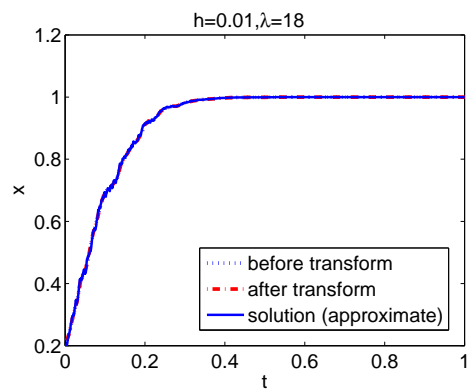}}
\\
\subfigure[$h=0.01,\quad\lambda=200$]{\label{f2.3}\includegraphics[width=0.5\textwidth]{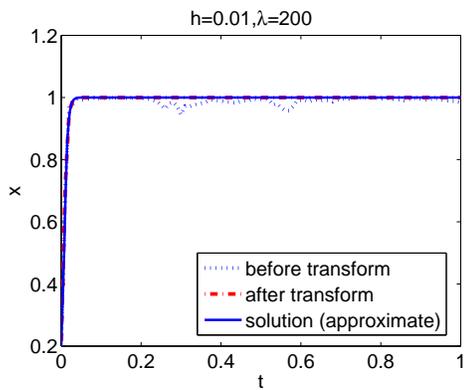}}
\hspace{1mm}
\subfigure[$h=0.001,\quad\lambda=200$]{\label{f2.4}\includegraphics[width=0.5\textwidth]{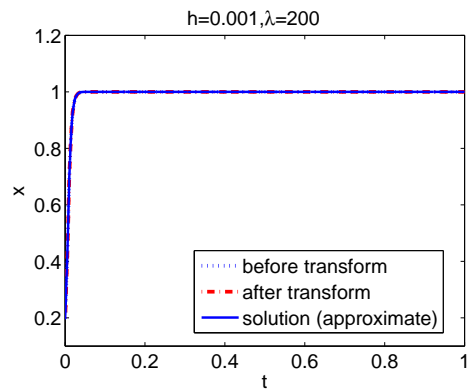}}
\caption{Comparison of stiffness before and after transformation}\label{f2}
\end{figure*}

In this section we compare the effect of numerical simulation of the stochastic protein kinetic equation based on the original equation (\ref{1.1}) and its transformed equation (\ref{2.4}), respectively. We apply the stochastic midpoint rule to both equations, with varying time step sizes for $\lambda=18$ and $\lambda=200$. The results are shown in Fig. \ref{f2}, for which we take $x_0=0.2$, $\sigma=1$, and the number of iterations in each time step for the realization of the implicit stochastic midpoint rule is 10.

In Fig. \ref{f2.1} and \ref{f2.2}, $\lambda=18$. As $h=0.12$, the midpoint rule applied to the equation before transform, i.e., the original equation (\ref{1.1}) loses accuracy severely, as shown by the blue dotted line in Fig. \ref{f2.1}, while that applied to the equation after transform, i.e., the equation (\ref{2.4}) together with the inverse transformation (\ref{2.6}) produces much better simulation as illustrated by the red dash-dotted line in the same figure. However, as $h$ is small enough, e.g. $h=0.01$, the midpoint rule applied to both equations gives good numerical results, which can be seen in Fig. \ref{f2.2}.

As the stiffness increases by enlarging the absolute value of $\lambda$, for example, $\lambda=200$, we see in Fig. \ref{f2.3} that the time step size $h=0.01$ loses its effectiveness if the numerical simulation is performed on the original equation (\ref{1.1}) without transformation, while remains valid if the numerical simulation is performed on the transformed equation (\ref{2.4}). However, for a much smaller time step size $h=0.001$, both simulations work fairly well again, with the cost of much more computations, as shown by Fig. \ref{f2.4}.

Note that, there is no explicit true solution for (\ref{1.1}) or (\ref{2.4}). The solution lines (blue solid) are simulated by the midpoint rule approximation based on the original equation (\ref{1.1}) with $h=0.001$ in Fig. \ref{f2.1}- \ref{f2.3}, and $h=0.0001$ in Fig. \ref{f2.4}.

To conclude, the numerical results show superiority of reducing stiffness of the SDE (\ref{1.1}) by variable transformation in the numerical computation of the equation. It permits larger time step sizes, and therefore reduces computational costs and increases computational robustness.
\section{Discussion}
\begin{figure*}
\subfigure[$\alpha=0.3,\quad h=0.12,\quad\lambda=18$]{\label{f3.1}\includegraphics[width=0.5\textwidth]{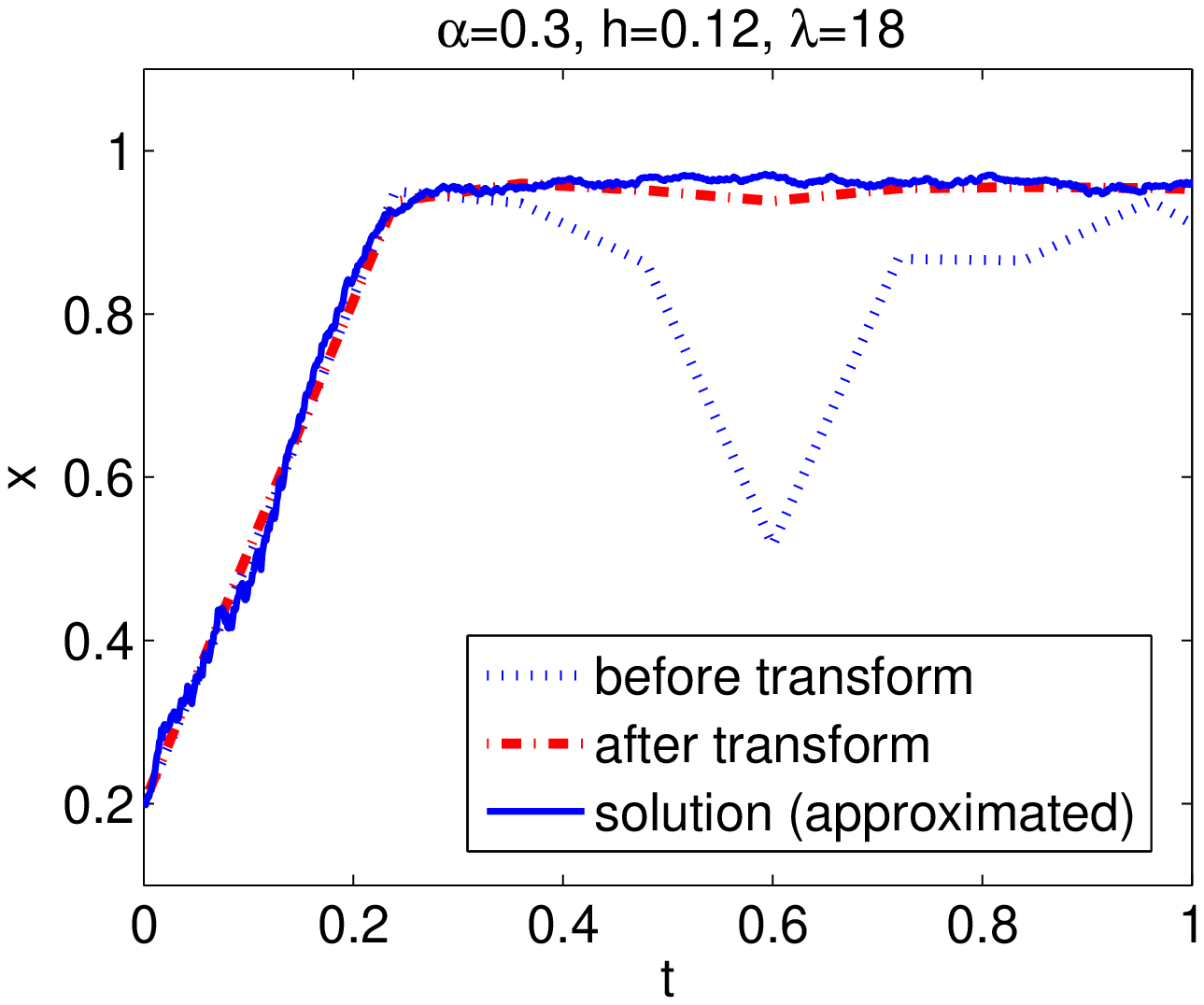}}
\hspace{1mm}
\subfigure[$\alpha=0.7,\quad h=0.12,\quad\lambda=18$]{\label{f3.2}\includegraphics[width=0.5\textwidth]{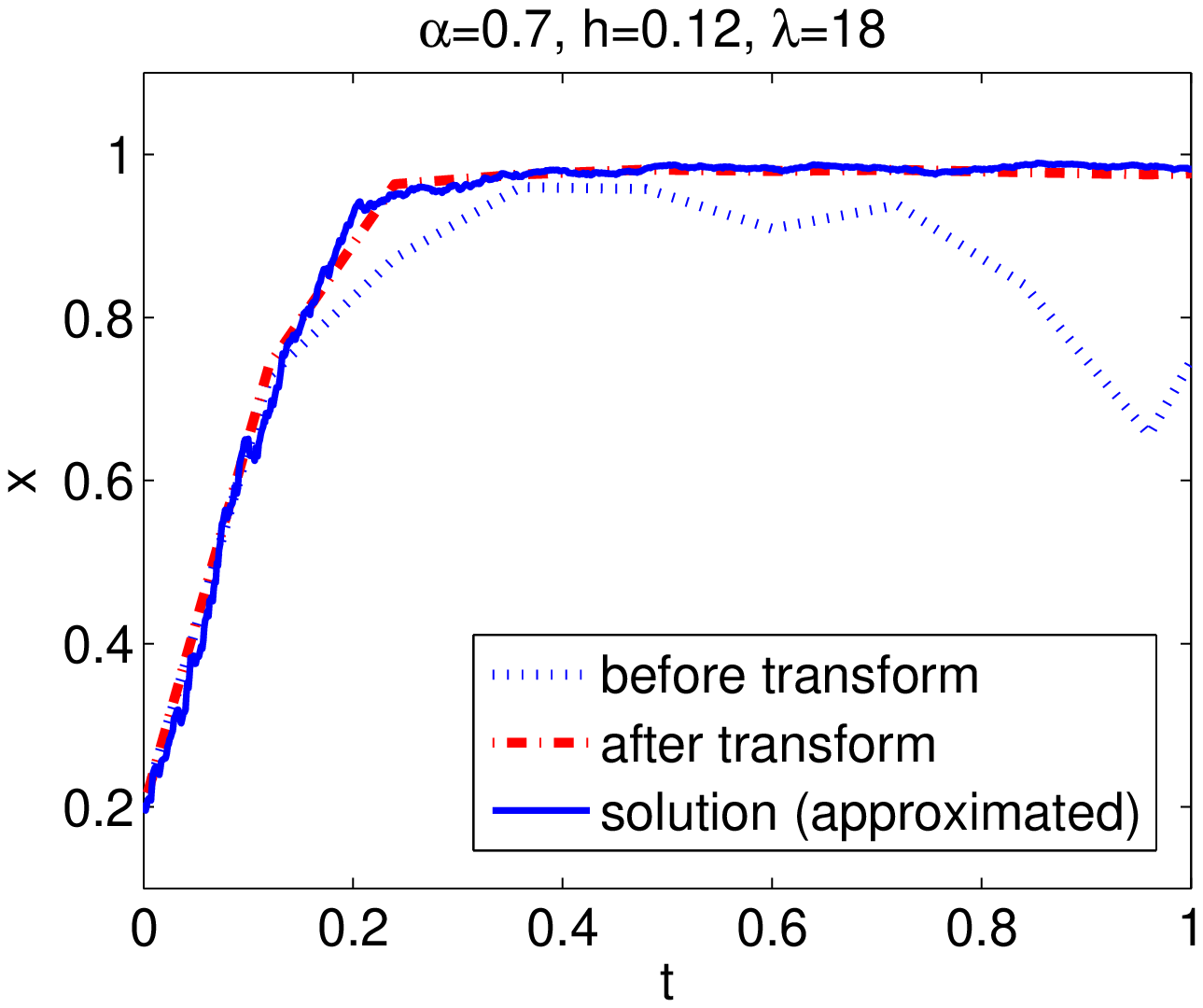}}
\caption{Stiffness reduction via transformation for $\alpha\in (0,1)$}\label{f3}
\end{figure*}

The more general stochastic protein kinetic equations involve a varying parameter $\alpha>0$ in the drift part
\begin{equation}
dx=(\alpha-x+\lambda x(1-x))dt+\sigma x(1-x)\circ dW(t),\quad x(0)=x_0,
\end{equation}
in which case $x=1$ is not a stochastic stationary solution if $\alpha\neq1$. However, if we still employ the variable transformation (\ref{2.5}), obtaining the transformed equation
\begin{equation}
\begin{split}
dX&=(\frac{1-\alpha}{1+\lambda}X^{-\lambda}-X+\frac{\lambda}{1+\lambda}X^{2+\lambda})dt\\&-(\frac{\sigma}{1+\lambda}X-\frac{\sigma}{1+\lambda}X^{2+\lambda})\circ dW(t)
\end{split}
\end{equation}
instead of (\ref{2.4}), we can still observe the effect of stiffness-reduction in the numerical tests for $\alpha\in(0,1]$, as shown by Fig. \ref{f3}.

As $|\lambda|$ is small, or $\lambda\ll -1$, the transformation is not recommended. In fact, this method of reducing stiffness can be generalized to a more general class of SDEs, e.g., of the form
\begin{equation}\label{4.1}
dy=f(y)dt+g(y)\circ dW(t), \quad y(0)=y_0,
\end{equation}
where the functions $f(y)$ and $g(y)$ are regular enough for the following discussion. This class of equations may cover many stochastic differential equation models in biology, chemistry, physics and so on.

Assume that there exists a constant $c$ such that $f(c)=g(c)\equiv 0$, and $y=c$ is a stochastically stationary solution of (\ref{4.1}). Thus, with $z=y-c$, the linearized equation of (\ref{4.1}) is (\cite{kandp})
\begin{equation}\label{4.2}
dz=f'(c)zdt+g'(c)z\circ dW(t).
\end{equation}
If $|f'(c)|$ is large, then the equation (\ref{4.2}) is stiff, and so is the original equation (\ref{4.1}). We search for a variable transformation $Z=Z(z)$ such that
\begin{equation}\label{4.5}
\frac{dZ}{dz}f'(c)z=-Z,
\end{equation}
which is an ordinary differential equation with solution
$|Z|=C|z|^{-\frac{1}{f'(c)}}, \quad C>0.$
For the purpose of having inverse transformation, we take
$Z=z^{-\frac{1}{f'(c)}}$
for $z>0$, and
$Z=(-z)^{-\frac{1}{f'(c)}}$
for $z<0$. In both cases we have the following transformed equation of (\ref{4.2})
\begin{equation}\label{4.9}
dZ=-Zdt-\frac{g'(c)}{f'(c)}Z\circ dW(t),
\end{equation}
which is less stiff than (\ref{4.2}) for large $|f'(c)|$. Based on these, we perform the transformation
$Y=(y-c)^{-\frac{1}{f'(c)}}$
for $y>c$ and
$Y=(c-y)^{-\frac{1}{f'(c)}}$
for $y<c$ on the equation (\ref{4.1}), to get its transformed equation
\begin{equation}\label{4.12}
\begin{split}
dY&=-\frac{1}{f'(c)}Y^{1+f'(c)}f(c+Y^{-f'(c)})dt\\&-\frac{1}{f'(c)}Y^{1+f'(c)}g(c+Y^{-f'(c)})\circ dW(t)
\end{split}
\end{equation}
for $y>c$ and
\begin{equation}\label{4.13}
\begin{split}
dY&=\frac{1}{f'(c)}Y^{1+f'(c)}f(c-Y^{-f'(c)})dt\\&+\frac{1}{f'(c)}Y^{1+f'(c)}g(c-Y^{-f'(c)})\circ dW(t)
\end{split}
\end{equation}
for $y<c$. Note that, if $y-c$ switches sign in the time interval of observation, then the transformation has no global inverse, which is a limitation of this method.
To show that (\ref{4.12}) or (\ref{4.13}) is less stiff than the original equation (\ref{4.1}), we need to have the linearized equation of them at $Y=0$, which is just the equation (\ref{4.9}) for both (\ref{4.12}) and (\ref{4.13}), for $f'(c)<0$. The less stiffness of (\ref{4.9}) than (\ref{4.2}) implies the less stiffness of (\ref{4.12}) and (\ref{4.13}) than (\ref{4.1}).

Note that the condition $f'(c)<0$ is also to guarantee that the drift parts of both (\ref{4.12}) and (\ref{4.13}) tends to $0$ as $Y\rightarrow 0$, and that $y=c$ is a stochastically asymptotically stable solution of (\ref{4.1}).
\section*{Acknowledgments}
The author is supported by the NNSFC (No.11071251, No.91130003, No. 11471310) and the 2013 Headmaster Funds of UCAS.

\end{document}